\documentclass[11pt,twoside]{article}

\usepackage{fullpage}
\usepackage{epsf}
\usepackage{fancyhdr}
\usepackage{graphics}
\usepackage{graphicx}
\usepackage{psfrag}
\usepackage{caption}
\usepackage{subcaption}
\usepackage{color}
\usepackage{amssymb}
\usepackage{amsthm}
\usepackage{amsfonts}
\usepackage{amsmath}
\usepackage{amssymb}
\usepackage{hyperref}
\usepackage{algorithm}
\usepackage{algorithmic}
\usepackage{bbm}
\usepackage{tikz}
\usepackage{caption}
\usepackage{enumitem}

\usepackage{mathrsfs}

\usepackage{tikz,lipsum,lmodern}
\usepackage[most]{tcolorbox}

\usepackage{macros-arXiv}

\theoremstyle{definition}

\newlength{\widebarargwidth}
\newlength{\widebarargheight}
\newlength{\widebarargdepth}
\DeclareRobustCommand{\widebar}[1]{%
  \settowidth{\widebarargwidth}{\ensuremath{#1}}%
  \settoheight{\widebarargheight}{\ensuremath{#1}}%
  \settodepth{\widebarargdepth}{\ensuremath{#1}}%
  \addtolength{\widebarargwidth}{-0.3\widebarargheight}%
  \addtolength{\widebarargwidth}{-0.3\widebarargdepth}%
  \makebox[0pt][l]{\hspace{0.3\widebarargheight}%
    \hspace{0.3\widebarargdepth}%
    \addtolength{\widebarargheight}{0.3ex}%
    \rule[\widebarargheight]{0.95\widebarargwidth}{0.1ex}}%
  {#1}}

\makeatletter
\long\def\@makecaption#1#2{
        \vskip 0.8ex
        \setbox\@tempboxa\hbox{\small {\bf #1:} #2}
        \parindent 1.5em  %% How can we use the global value of this???
        \dimen0=\hsize
        \advance\dimen0 by -3em
        \ifdim \wd\@tempboxa >\dimen0
                \hbox to \hsize{
                        \parindent 0em
                        \hfil
                        \parbox{\dimen0}{\def\baselinestretch{0.96}\small
                                {\bf #1.} #2
                                }
                        \hfil}
        \else \hbox to \hsize{\hfil \box\@tempboxa \hfil}
        \fi
        }
\makeatother

\long\def\comment#1{}

\newcommand{\red}[1]{\textcolor{red}{#1}}

\newcommand{\Exs}{\ensuremath{{\mathbb{E}}}}

\newcommand{\grad}[1]{\ensuremath{\nabla #1}}

 \newcommand{\f}{\ensuremath{f}}

\newcommand{\real}{\ensuremath{\mathbb{R}}}

\newcommand{\xstar}{\ensuremath{x^{*}}}

\newcommand{\TotalEpochs}{\ensuremath{M}}

\newcommand{\nobs}{\ensuremath{N}}

\newcommand{\stepsize}{\ensuremath{\lambda}}
\newcommand{\epochLen}{\ensuremath{{K}}}

\newenvironment{carlist}
 {\begin{list}{$\bullet$}
 {\setlength{\topsep}{0in} \setlength{\partopsep}{0in}
  \setlength{\parsep}{0in} \setlength{\itemsep}{\parskip}
  \setlength{\leftmargin}{0.15in} \setlength{\rightmargin}{0.08in}
  \setlength{\listparindent}{0in} \setlength{\labelwidth}{0.08in}
  \setlength{\labelsep}{0.1in} \setlength{\itemindent}{0in}}}
 {\end{list}}

 \newcommand{\bcar}{\begin{carlist}}
\newcommand{\ecar}{\end{carlist}}

\usepackage{array}
\newcolumntype{P}[1]{>{\centering\arraybackslash}p{#1}}

\newcommand{\Ncal}{\mathcal{N}}

\newcommand{\NumEpoch}{\ensuremath{\TotalEpochs}}
\newcommand{\epiter}{\ensuremath{m}}

\newcommand{\xhat}{\widehat{x}}
\newcommand{\epoch}{m}
\newcommand{\xbar}{\widebar{x}}
\newcommand{\epochsize}{T}
\newcommand{\vargrad}{\widetilde{\nabla} \fun}

\newcommand{\gradfunbar}{\widebar{\nabla} \f} 
\newcommand{\ftilde}{\tilde{\f}}
\newcommand{\gtilde}{\tilde{g}}

\newcommand{\algo}{\mathcal{A}}

\newcommand{\projT}{\mathrm{P}_{\mathcal{T}}}
\newcommand{\Hess}{H^\star}
\newcommand{\CovMat}{\Sigma^\star}
\newcommand{\Cov}{\mathrm{Cov}}

\newcommand{\pen}{\mathrm{R}}
\newcommand{\fun}{\f}

\newcommand{\Problem}{\mathcal{P}}
\newcommand{\Xset}{\mathcal{X}}
\newcommand{\Pdist}{\ensuremath{\mathbb{P}}}
\newcommand{\muPar}{\mu}
\newcommand{\Lipcon}{\mathrm{L}}

\newcommand{\g}{g}
\newcommand{\dims}{d}
\renewcommand{\xstar}{x_{\P_0}^\star}

\begin{document}

%%%%%%% Title PAGE %%%%%%%%%%%%%%%%%%%%%%%%%%%%%%%%%%%%%%%%%%%%%%%%%%%

\begin{center}

{\bf{\LARGE{\mbox{Stochastic Optimization with Constraints:} \\[.1cm]
      A Non-asymptotic Instance-Dependent Analysis}}}

\vspace*{.2in}

{\large{
 \begin{tabular}{ccc}
  Koulik Khamaru$^{\diamond}$
 \end{tabular}
}}

 \vspace*{.2in}

 \begin{tabular}{c}
 Department of Statistics$^\diamond$ \\
 Rutgers University\\
 \end{tabular}

\vspace*{.2in}

\begin{abstract}
We consider the problem of stochastic convex optimization under convex constraints. We analyze the behavior of a natural variance reduced proximal gradient (VRPG) algorithm for this problem. Our main result is a non-asymptotic
guarantee for VRPG algorithm. Contrary to minimax worst case guarantees, our result is instance-dependent in nature. This means that our guarantee captures the complexity of the loss function, the variability of the noise, and the geometry of the constraint set.
We show that the non-asymptotic performance of the VRPG algorithm is governed by the scaled distance (scaled by $\sqrt{\nobs}$) between the solutions of the given problem and that of a certain small perturbation of the given problem --- both solved under the given convex constraints; here, $\nobs$ denotes the number of samples. 
Leveraging a well-established connection between local minimax lower bounds and solutions to perturbed problems, we show that as $\nobs \rightarrow \infty$, the VRPG algorithm achieves the renowned local minimax lower bound by  H\`{a}jek and {Le Cam} up to universal constants and a logarithmic factor of the sample size.
  
\end{abstract}

\end{center}

\section{Introduction}
\label{sec:intro}
In this paper,  we study the problem of stochastic convex optimization problem under constraints. Concretely, we are interested in problems of the form 
\begin{align}
\label{eqn:constraint-opt}
\min_{x \in \Xset} \;  \f(x) &\equiv  \Exs_{z \sim \Pdist_0} \f(x, z)  \notag \\
\qquad   \text{where} \;\; \Xset &= \left\lbrace x \in \real^\dims \; : \; \g_i(x) \leq 0 \;\;
\text{for} \;\; i = 1, \ldots, m \right\rbrace.  
\end{align}
Here, $\g_i$ are convex and twice continuously differentiable functions, and the function $\f$ is twice continuously differentiable, $\mu$-strongly convex and $\Lipcon$-smooth for some $0 < \mu \leq \Lipcon$. Given $\nobs$ iid samples $z_1, z_2, \ldots, z_\nobs \sim \P_0$, we are interested in designing computationally efficient algorithm that outputs an approximate solution of the problem~\eqref{eqn:constraint-opt}.

Constraint stochastic optimization problems of the form~\eqref{eqn:constraint-opt} are ubiquitous in statistics and optimization, and the study of these problems dates backs to more than 70 years. 
There are two standard strategies for solving the problem~\eqref{eqn:constraint-opt}: (a) M-estimators under constraints (b) Stochastic approximation methods.  

\subsection{M-estimators under constraints:}
The first approach, popular in the statistics community, provides an approximate solution of problem~\eqref{eqn:constraint-opt} by solving a constrained M-estimation problem. Given $\nobs$ iid samples $z_1, z_2, \ldots, z_\nobs \sim \P_0$, one outputs the solution of the following M-estimation problem
\begin{align}
\label{eqn:M-estimator}
\min_{x \in \Xset} \frac{1}{\nobs} \sum_{i = 1}^\nobs \f(x, z_i).
\end{align}  
The asymptotic behavior of the above problem is well understood by now. The study of this problem has its roots in the works of Dupacoka and Wets~\cite{dupacova1988asymptotic}, Shapiro~\cite{shapiro1989asymptotic}, and  King and Rockafellar~\cite{king1993asymptotic}. In a recent work~Duchi and Ruan~\cite{duchi2021asymptotic} and Davis et al.~\cite{davis2023asymptotic} show that the solution of the solution obtained from~\eqref{eqn:M-estimator} is asymptotically locally minimax optimal in the sense of H\'{a}jek and Le Cam~\cite{le1972limits,hajek1972local,le2000asymptotics,van2000asymptotic,wellner2013weak}. In simple words,
in a local shrinking neighborhood (shrinking as the number of the sample size increases) of the problem,
 there is no other algorithm that performs better than the solution of the problem~\eqref{eqn:M-estimator}. Contrary to global minimax bounds which capture the difficulty of the hardest problem in the problem class, the local minimax lower bounds from~\cite{duchi2021asymptotic} captures a more finer notion of difficulty of the specific problem at hand. 

\subsection{Stochastic approximation methods:}       
While the asymptotic performance of constrained M-estimator problems are optimal asymptotically, one drawback of such methods is that they may be expensive to compute. This motivated people to design computationally attractive stochastic approximation schemes that converge to the solution of problem~\eqref{eqn:constraint-opt}.

\newcommand{\myc}[1]{\textcolor{red}{#1}}

\paragraph{Suboptimality of dual-averaging:} The problem of  stochastic optimization under constraints has been has been studied in the optimization community for some time; see the work of  and  King~\cite{king1986asymptotic}. 
The behavior stochastic algorithms are more nuanced due to presence of constraint set $\Xset$. For instance, Duchi and Ruan~\cite{duchi2021asymptotic} show that stochastic approximation with
dual averaging~\cite{nesterov2009primal} --- a natural algorithm in this case --- is not always optimal. Concretely, in~Section 5.2 of their paper Duchi and Ruan argue that the dual averaging algorithm is optimal when the the constraint set $\Xset$ consists of linear constraints, but the same algorithm is sub-optimal when $\Xset$ contains non-linear constraints.   

\paragraph{Asymptotic optimality via Polyak-Ruppert averaging:}
In a surprising new result, Davis, Drusvyatskiy, and Jiang~\cite{davis2023asymptotic} show that vanilla stochastic approximation algorithm when combined with a Polyak-Ruppert (PR) averaging step is asymptotically local-minimax  optimal. Staring with an $x_1 \in \mathcal{X}$, the updates of the PR-averaged stochastic approximation as obtained as follows: 
\begin{align*}
\texttt{(SA-step)} \qquad x_{k+1} &= \text{Proj}_{\mathcal{X}}\left(x_k - \alpha_k \cdot \grad \f(x_{k}, z_k) \right). \\
\texttt{(PR-averaging)} \qquad \widebar{x}_{k + 1}& = \frac{1}{k + 1} \sum_{\ell = 1}^{k + 1} x_\ell 
\end{align*} 
Davis et al.~\cite{davis2023asymptotic} argue that under some mild assumptions of the loss function $\f$ and the constraint set $\mathcal{X}$, the scaled difference $\sqrt{\nobs}(\widebar{x}_\nobs - \xstar)$ converges to a normal distribution with the optimal covariance matrix. Here, $\xstar$ is the (unique) solution to the problem~\eqref{eqn:constraint-opt}, and $\nobs$ is the number of samples used (or the number of steps in the algorithm). The authors also argue that $\sqrt{\nobs}\Exs \|\widebar{x}_\nobs - \xstar\|_2$ matches a local minimax lower bound suggested by the H\'ajek Le Cam theory~\cite{hajek1972local,le1972limits,duchi2021asymptotic}.   

\subsection{Non-asymptotic bounds}
While the asymptotic performances of both M-estimator and stochastic approximation algorithms are quite well understood, not much is known
about the non-asymptotic instance-dependent performance of either methods for the problem~\eqref{eqn:constraint-opt} when constraint set $\mathcal{X} \neq \real^\dims$. In a recent work~\cite{duchi2016local}, the authors provide a non-asymptotic lower bounds for the problem~\eqref{eqn:constraint-opt} under some growth condition on an appropriate modulus of continuity. In particular, motivated by the work of~\cite{cai2015framework}, the authors in the work~~\cite{duchi2016local} construct a non-asymptotic local minimax lower bound for problem~\eqref{eqn:constraint-opt} by comparing problem~\eqref{eqn:constraint-opt} with its hardest alternative in a local shrinking  neighborhood (shrinking with sample size $\nobs$) of the problem. The authors also argue that in one dimension, a grid search based algorithm attains this non-asymptotic local lower bound. Unfortunately, in higher dimensions, grid search algorithms become computationally prohibitive. Thus it is not clear whether such results can be extended to higher dimension.

\subsection{Contribution and organization}
Motivated from the recent success of instance dependent understanding of stochastic approximation algorithms, our goal in this paper is to have a non-asymptotic understanding of stochastic approximation algorithms in a instance-dependent sense. Put simply, we would like to derive non-asymptotic upper-bounds that reflect the structure of the loss function $\f$, the variability of the noise, and the structure of the constraint set $\mathcal{X}$. In Section~\ref{sec:benchmark},  we motivate non-asymptotic instance dependence benchmark for the problem~\ref{eqn:constraint-opt}. This benchmark is motivated from the well-studied connection between the local minimax lower bounds and solutions of perturbed problem~\cite{duchi2021asymptotic,davis2023asymptotic}.  Section~\ref{sec:pen-class} is devoted to a devoted to a detailed description of the problem that we consider. 
In Section~\ref{SecMain} we state the variance reduced proximal gradient (VRPG) algorithm that we study (see Algorithm~\ref{Algo:VR-prox}), and in Theorem~\ref{thm:VRPG}, we provide non-asymptotic instance-dependent upper bound for
the VRPG algorithm. Finally, we argue that  VRPG is asymptotically locally minimax optimal up to a factor which is logarithmic in the number of sample size $\nobs$ and universal constants.

\section{An instance-dependent local minimax lower bound}
\label{sec:benchmark}
We start with a discussion on the statistical limit of the constrained stochastic optimization problem~\eqref{eqn:constraint-opt}. Our lower bounds are motivated from the local minimax theory of Le Cam and H\'{a}jek~\cite{le1972limits,hajek1972local}. Contrary to ,standard minimax guarantees --- which capture the difficulty of the worst case problem in a problem class, local minimax guarantees captures a problem specific difficulty of the problems. In other words, these lower bounds are adaptive to the structure of the function $\f$ as well as the constrained set $\Xset$ in problem~\eqref{eqn:constraint-opt}. 

The idea is to consider a small neighborhood of the given problem $\Problem_0$ --- shrinking with the sample size --- and consider worst case difficulty in this local neighborhood. Intuitively, in a small local neighborhood of a given problem $\Problem_0$, problem difficulty does not change drastically. Consequently, the worst case difficulty in the local neighborhood of the problem $\Problem_0$ is a finer notion of difficulty of problem $\Problem_0$ compared to the worst case problem difficulty.  
\newcommand{\KL}{\mathrm{KL}}

Focusing on our problem~\eqref{eqn:constraint-opt}, let $\P$ is any other distribution over the randomness $z$ such that $\KL(\P, \P_0)$ is small to ensure 
\begin{align}
\f_{\Pdist} = \Exs_{z \in \Pdist} \f(x,z) \quad \text{is strongly convex}. 
\end{align}
Let $x_{\P}^\star$ the unique solution to the perturbed problem: 
\begin{align}
\label{eqn:constraint-opt-tilted}
x_{\P}^\star = \arg\min_{x \in \Xset} \; \left\lbrace \f_{\Pdist}(x) \equiv  \Exs_{z \sim \Pdist} \f(x, z) \right\rbrace. 
\qquad   
\end{align}
Given access to $\nobs$ iid samples $z_1, \ldots, z_\nobs$ from the distribution $\Pdist$, the following result from the paper \cite[Corollary 1]{duchi2021asymptotic} characterizes the local minimax lower bound for the problem~\eqref{eqn:constraint-opt}. 

\vspace{10pt}

\begin{tcolorbox}
\begin{proposition}[Informal]
\label{prop:local-minimax-lb}
Let $\xhat_\nobs$ be the output of \emph{any} algorithm with access to $\nobs$ iid samples $z_1, \ldots, z_\nobs$. Then under suitable smoothness condition on the function $\f$ and the constraint functions $\g_i's$ in problem~\eqref{eqn:constraint-opt}, we have the following asymptotic lower bound:
\begin{align}
\label{eqn:minimax-lb}
\liminf_{c \rightarrow \infty}  \liminf_{\nobs \rightarrow \infty} \sup_{\Pdist: \KL(\Pdist || \Pdist_0) \leq c/\nobs} 
\Exs_{\Pdist^\nobs} \left[\|\sqrt{\nobs}(\xhat_\nobs - x_{\P}^\star)\|_2^2\right] \geq \Exs\|Z\|_2^2
\end{align}
where $Z \sim \Ncal(0, \projT \Hess \projT \CovMat \projT \Hess \projT)$. 
\end{proposition}
\end{tcolorbox}

\newcommand{\lagrange}{\mathcal{L}}

\vspace{10pt}
\noindent 
A few comments regarding the last lower bound are in order. The lower bound depends on the property of three matrices $\projT, \Hess$ and $\CovMat$. Let us understand these three matrices one at a time. We use $\lagrange(\cdot, \cdot)$ as the Lagrangian of the problem~\eqref{eqn:constraint-opt}, i.e., 
\begin{align}
\label{eqn:lagrangian}
\lagrange(x, \beta) = \f(x) + \sum_{i = 1}^{m} \beta_i \g_i(x).
\end{align}

\begin{enumerate}
\item[(a)] 
The matrix $\Hess$ in the bound~\eqref{eqn:minimax-lb} is the hessian of the Hessian of the Lagrangian~\eqref{eqn:lagrangian} at optima $\xstar$
\begin{align}
\Hess = \grad^2 \lagrange(\xstar, \beta^\star) = 
\grad^2 \f(\xstar) + \sum_{i = 1}^m \beta_i^\star \grad^2 \g_i(\xstar)
\end{align}
where $\beta_i^\star \geq 0$ with $\beta_i^\star = 0$ when $\grad \g_i(\xstar) < 0$ are unique scalars satisfying 
\begin{align*}
\grad \f(\xstar) + \sum_{i = 1}^m \beta_i^\star \g_(\xstar) = 0
\end{align*}
See the paper~\cite[Section 2]{duchi2021asymptotic}, and the book~\cite[Chapter VII.2]{hiriart1996convex} for more a more detailed discussion on existence of such $\beta^\star$. 
Put simply, the Hessian matrix $\Hess$ captures the structures of the loss function $\f$ and the constraint-set $\Xset$ near optima $\xstar$. 
\item[(b)] 
The matrix 
$\CovMat = \Cov_{z \in \Pdist_0}(\grad \f(x_{\Pdist_0}^\star, z))$ is the covariance matrix of the function gradient evaluated at the optima $x_{\Pdist_0}^\star$. This term captures the noise present in the problem~\eqref{eqn:constraint-opt}. 
\newcommand{\Cone}{\mathcal{T}}
\item[(c)] Let $\Cone_{\Xset}(\xstar)$ denote the critical tangent cone at $\xstar$: 
\begin{align*}
\Cone_{\Xset}(\xstar) &:= \left\lbrace w \in \real^\dims \; : \; 
\grad \g_i(\xstar)^\top w = 0 \;\; \text{for all} \;\; i \;\; \text{with} \;\; \g_i(\xstar) = 0  \right\rbrace \\
&= \cap_{i = 1}^{m_0} \left\lbrace w \in \real^\dims \; : \;  w \in \real^\dims \; : \; 
\grad \g_i(\xstar)^\top w = 0  \right\rbrace     
\end{align*} 
Here, we assumed, without loss of generality, that the first $m_0$ constraints $\left\lbrace \g_i(x) \leq 0\right\rbrace$ are active at the optima $\xstar$. The matrix $\projT$ denote the projection matrix to this tangent cone $\Cone_{\Xset}(\xstar)$. This term captures the structure of the constraint set $\Xset$ near the optima $\xstar$.
\end{enumerate}
The lower bound from Proposition~\ref{prop:local-minimax-lb} provides an asymptotic lower bound for any procedure, and it is not immediately clear what could be a non-asymptotic benchmark. To motivate our non-asymptotic benchmark we discuss an alternative description of the lower bound in Proposition~\ref{prop:local-minimax-lb}. We do so by connecting the lower bound~\eqref{eqn:minimax-lb} with solutions of a class of perturbed problems obtained from the original problem~\eqref{eqn:minimax-lb}.

\subsection{Perturbed problems and connections to local minimax lower bounds}
Before we discuss the connection between perturbations of the problem~\eqref{eqn:constraint-opt} and the local minimax lower bound~\eqref{eqn:minimax-lb}, it is helpful to first understand which distributions $\Pdist$ with $\KL(\Pdist || \Pdist_0)$ achieve the lower bound in Proposition~\ref{prop:local-minimax-lb}. It turns out that the lower bound is achieved by considering the \emph{tilted} distribution: 
\begin{align}
\label{eqn:tilted-dist}
d\Pdist_u(z) &= \frac{1 + \left\langle u, \grad \f(x_{\P_0}^\star,z) - \grad \f(x_{\P_0}^\star) \right\rangle}{C_0} \cdot d\Pdist_0(z) \quad \text{where} \notag \\
 \quad C_0 &= 1 + 
\int \left\langle u, \grad \f(x_{\P_0}^\star,z) - \grad \f(x_{\P_0}^\star) \right\rangle \cdot d \P_0(z).
\end{align}
and taking supremum over the set $\|u\|_2 \leq c/\sqrt{\nobs}$; (see Theorem~1 in the paper~\cite{duchi2021asymptotic} for additional details).
Put simply, 
\begin{center}
\emph{The hardest alternative to distribution $\P_0$ is obtained when we \emph{tilt} the distribution $\P_0$ in the direction $[\grad \f(x_{\P_0}^\star,z) - \grad \f(x_{\P_0}^\star)]$ --- the noise present in the sample function gradient}.
\end{center}
 Any distribution $\P_u$ of the form~\eqref{eqn:tilted-dist} produces a function $\f_{u} := \Exs_{z \in \P_u} \f(x, z)$ which is a linear perturbation to the original function $\f$. 
The study of the change in the optimal solution of the optimization problem~\eqref{eqn:constraint-opt} and its \emph{tilted} counterpart $\f_u$ is related to the study of tilt-stability~\cite{poliquin1998tilt} in variational analysis.  The idea of tilt-stability has also been very influential in developing asymptotic results for the problem~\eqref{eqn:constraint-opt}; see for instance the works by Shapiro~\cite{shapiro1989asymptotic}. These ideas are also building blocks for deriving local asymptotic lower bounds~\cite{duchi2021asymptotic,davis2023asymptotic}.

\vspace{10pt}
\noindent 
\subsubsection{A particular perturbation of interest}
A particular perturbation of interest is the following. 
Given $\nobs$ iid samples $z_1, \ldots, z_\nobs$, we define a tilted function $\f_{\nobs}$ as   

\begin{gather}
\f_{\nobs}(x) = \f(x) + \left\langle x, \widebar{\nabla} \f(x_{\P_0}^\star) - \grad \f(x_{\P_0}^\star) \right\rangle \quad \text{where}
\quad 
\widebar{\nabla} \f(x_{\P_0}^\star)
 = \frac{1}{\nobs} \sum_{i = 1}^\nobs \nabla \f(x_{\P_0}^\star, z_i)
 \notag  \\
  x_\nobs^\star = \arg\max_{x \in \Xset} \f_{\nobs}(x) 
  \label{eqn:tilded-EMP-loss}.  
\end{gather}

\noindent This particular perturbation~\eqref{eqn:tilded-EMP-loss} is motivated from the behavior of the classical M-estimator. Indeed, let 
\begin{align*}
\text{(M-estimator:)} \qquad \qquad  \widehat{x}^{M} := \arg\min_{x \in \Xset} \left\lbrace  \widehat{\f}_\nobs(x):= \frac{1}{\nobs} \sum_{i = 1}^\nobs \f(x ; z_i)  \right\rbrace, \qquad \qquad \qquad\qquad\qquad 
\end{align*}
Then under some smoothness assumption on the function $\f$, we have
\begin{align*}
\grad \widehat{\f}_\nobs(x) =  \grad \f_\nobs(x) + o_{p}(1) \cdot \left(x - \xstar \right)  
\end{align*}
The above relation ensures that the the asymptotic behavior of the M-estimator $\widehat{x}^{M}$ is same as that of $x_\nobs^\star$ from~\eqref{eqn:tilded-EMP-loss}. See Section 2.3 from the paper~\cite{duchi2021asymptotic} and also the paper~\cite{shapiro1989asymptotic} for details. The perturbed problem $\f_\nobs$ from~\eqref{eqn:tilded-EMP-loss} lies at the heart of variance reduced stochastic optimization algorithms~\cite{johnson2013accelerating,wainwright2019variance,khamaru2021temporal,khamaru2021instance}. We are now ready to connect the perturbed problem~\eqref{eqn:tilded-EMP-loss} with the local minimax lower bound~\eqref{eqn:minimax-lb}, and this connection forms the basis of the non-asymptotic benchmark that we consider. 

Intuitively, as $\nobs$ increases, the function $\f_\nobs$ becomes increasingly close to $\f$, and consequently, $x_\nobs^\star$ becomes close to $x_{\P_0}^\star$. Proposition~1 from the paper~\cite{duchi2021asymptotic} formalizes this intuition:
\begin{align*}
x_\nobs^\star = x_{\P_0}^\star + 
\projT \Hess \projT  v + o(\|v\|)
\quad \text{where} \quad v = \widebar{\nabla} \f(x_{\P_0}^\star) - \grad \f(x_{\P_0}^\star).   
\end{align*}
Also, see the book by~\cite{dontchev2009implicit} for a detailed discussion. Now, scaling the difference between $x_\nobs^\star$ and $x_{\P_0}^\star$ by $\sqrt{\nobs}$ we obtain the following result
\vspace{10pt}
\begin{tcolorbox}
Under the assumptions of Proposition~\ref{prop:local-minimax-lb}
\begin{align}
\label{eqn:perturbation-bound}
\lim_{\nobs \rightarrow \infty} \sqrt{\nobs} (x_\nobs^\star - x_{\P_0}^\star)  \stackrel{d}{\rightarrow} \Ncal(0, \projT \Hess \projT \CovMat \projT \Hess \projT)  
\end{align}
Combining the statement~\eqref{eqn:perturbation-bound} and Proposition~\eqref{prop:local-minimax-lb} we observe that the local minimax lower bound for the problem~\eqref{eqn:constraint-opt} is characterized by 
$\sqrt{\nobs} (x_\nobs^\star - x_{\P_0}^\star)$ --- the (scaled) difference between the optimal solutions of the function $\f$ and its perturbed counterpart $\f_\nobs$. 
\end{tcolorbox}

\subsection{A non-asymptotic instance-dependent benchmark:}
\label{sec:benachmark}
Motivated from the discussion in the last two sections, we use $\delta_{\nobs}$ as the scaled difference between $x_\nobs^\star$ and $x_{\P_0}^\star$:
\begin{align}
\label{eqn:benchmark-1}
\delta^2(\nobs) := \nobs \cdot \Exs\|x_\nobs^\star - x_{\P_0}^\star\|_2^2
\end{align} 
Here, $x_\nobs^\star$ and $x_{\P_0}^\star$ are defined in equations~\eqref{eqn:tilded-EMP-loss} and~\eqref{eqn:constraint-opt-tilted} respectively. A few comments regarding the benchmark $\delta^2(\nobs)$ are in order. 

\begin{enumerate}
\item[(a)] The benchmark $\delta^2(\nobs)$ is non-asymptotic in nature as it depends on number of samples available.   
\item[(b)] The perturbed function $\f_\nobs$ depends on the 
randomness noise present in the problem. Consequently, the benchmark $\delta^2(\nobs)$ captures the noise in the problem. 
\item[(c)] The benchmark $\delta^2(\nobs)$ depends on the constraint-set $\Xset$ through $x_\nobs^\star$ and $x_{\P_0}^\star$. 
\end{enumerate}
In the rest of the paper, our goal is the following:
\vspace{10pt}
\begin{tcolorbox}
Can we design an algorithm $\algo$ such that given $\nobs$ samples, the algorithm outputs $\widehat{x}_{\algo}$ with  
\begin{align}
\label{eqn:UB-banchmark-aim}
\Exs \|\widehat{x}_{\algo} - x_{\P_0}^\star\|^2 \lesssim \frac{\delta^2(\nobs)}{\nobs} \qquad \text{for all} \;\; \nobs. 
\end{align} 
The symbol $\lesssim$ above bound indicate that the above inequality holds up to some universal constants and logarithmic factors. 
\end{tcolorbox}

\section{Problem setup}
\label{sec:pen-class}
While our discussions so far is motivated from the constraint optimization problem~\eqref{eqn:constraint-opt}, the ideas used to develop our algorithm also applies to a more general class of problems. In particular, we consider problems of the following form:
\begin{align}
\label{eqn:lagrange-problem}
\min_{x} \left\lbrace \Exs_{z \in \P_0} \f(x, z) + \pen(x) \right\rbrace 
\end{align}
Here, $\pen(x)$ is a known penalty function, which is assumed to be convex, proper and lower semi-continuous. The constraint optimization problem~\eqref{eqn:constraint-opt} is a special case of problem~\eqref{eqn:lagrange-problem} when we take $\pen(x) = 1_{x \in \Xset}$ --- the indicator function of the convex set $\Xset$. Throughout, we assume the following structural assumptions on the functions $\f(\cdot, \cdot)$ and $\pen$.

\subsection*{Assumptions}
\begin{enumerate}[label= \textbf{A\arabic*.}]
\item The loss-function $\fun(x) = \Exs_{z\sim \P_0}\f(x,z)$ is $\mu$-strongly convex and $\Lipcon$-smooth, for some $0 < \muPar \leq \Lipcon <\infty$.~\label{assn:mu-L-funcs}  
\item The penalty function $\pen$ is a  proper, convex and lower semi-continuous function. We assume that the function $\pen$ is known. ~\label{assn:penalty} 
\item The function $\f(\cdot, \cdot)$ is continuously differentiable in its first argument, and the gradient with respect to $x$, denoted as $\grad \f(x, z)$, is $\Lipcon$ Lipschitz continuous in its first argument, i.e., 
\begin{align*}
\|\grad \f(x, z) - \grad \f(y, z) \|_2 \leq \Lipcon \|x - y\|_2 \quad \text{for all} \;\; z.
\end{align*}
~\label{assn:Lipschitz-gradient} 
\end{enumerate}

\subsubsection*{Instance-dependent benchmark:}
Following the discussion of Section~\ref{sec:benchmark}, let $x_{\P_0}^\star$ be the \emph{unique} solution to the problem~\eqref{eqn:lagrange-problem}, i.e., 
\begin{align}
\label{eqn:lagrange-problem-2}
x_{\P_0}^\star =  \arg\min_{x} \left\lbrace \Exs_{z \in \P_0} \f(x, z) + \pen(x) \right\rbrace 
\end{align}
Given access to $\nobs$ iid samples $z_1,\ldots,z_\nobs$ from $\P_0$, we deine the perturbed function $\f_\nobs$ and solution $x_{\nobs}^\star$ as follows.
\begin{gather}
\f_{\nobs}(x) = \f(x) + \left\langle x, \widebar{\nabla} \f(x_{\P_0}^\star) - \grad \f(x_{\P_0}^\star) \right\rangle \quad \text{where}
\quad 
\widebar{\nabla} \f(x_{\P_0}^\star)
 = \frac{1}{\nobs} \sum_{i = 1}^\nobs \nabla \f(x_{\P_0}^\star, z_i)
 \notag  \\
  x_\nobs^\star = \arg\min_x  \left\lbrace \f_{\nobs}(x) + \pen(x) \right\rbrace \label{eqn:tilded-penalized-version}.  
\end{gather}
Note that the function $\f_{\nobs}$ is strongly convex by construction, and as a result the solution $x_\nobs^\star$ is unique. 
The benchmark $\delta^2(\nobs)$ is again defined as the scaled difference square-difference between $x_\nobs^\star$ and $\xstar$. 
\begin{align}
\label{eqn:benchmark-2}
\delta^2(\nobs) := \nobs \cdot \Exs\|x_\nobs^\star - x_{\P_0}^\star\|_2^2
\end{align}

\section{Main results}
\label{SecMain}
This section is devoted to the description of our algorithm and its performance guarantees. In Section~\ref{sec:VR-algo}, we first describe our variance reduced algorithm~\ref{Algo:VR-prox}, and discuss the intuition behind the algorithm. In Section~\ref{sec:VRPE-thm-section}, we provide non-asymptotic guarantees for the output of Algorithm~\ref{Algo:VR-prox}. 
\subsection{Variance-reduced projected stochastic gradient}
\label{sec:VR-algo}
Our proposed algorithm is a variance reduced stochastic projected gradient descent algorithm and is stated in Algorithm~\ref{Algo:VR-prox}. Just like standard variance reduction schemes~~\cite{johnson2013accelerating} our algorithm proceeds in epochs. 
\begin{enumerate}
\item[(a)] At the beginning of epoch, we create an approximation of the true constraint optimization problem~\eqref{eqn:lagrange-problem}.
\item[(b)] Within each epoch, we run a proximal stochastic gradient algorithm on the approximate problem constructed step (a). 
\end{enumerate}
Concretely, at epoch $\epoch$ we start with an estimate $\xbar_\epoch$ of the true minimizer $\xstar$ of the problem~\eqref{eqn:lagrange-problem}, and $\epochsize$ iid sample $\{z_1, \ldots, z_\epochsize\}$ samples from the distribution $\P_0$, we construct a modified estimate of gradient as follows 
\begin{align}
\label{eqn:recentered-problem}
\vargrad(x, z) := \nabla \f(x, z) + 
\left\lbrace \gradfunbar(\xbar_\epoch) - \grad \f (\xbar_\epoch, z) \right\rbrace
\quad \text{where} \quad \gradfunbar(\widebar{x}_\epoch) = \frac{1}{\epochsize} \sum_{i = 1}^\epochsize \grad \f(\widebar{x}_\epoch, z_i) 
\end{align}
where $z$ is independent of $\{z_1, \ldots, z_\epochsize\}$. By construction,
\begin{align}
\Exs[\vargrad(x, z)] = \grad \fun(x) \quad \text{for all} \quad x.   
\end{align}
Additionally, when $x$ is $\xbar_\epoch$ is \emph{close}, the variance of $\vargrad(x, z)$ is much smaller than $\grad \f(x, z)$. Indeed,
invoking the Lipschitz gradient assumption~\ref{assn:Lipschitz-gradient}, and using Cauchy-Schwartz inequality we have  
\begin{align*}
\Exs \| \vargrad(x, z) - \grad \fun (x) \|^2 &\leq 3 \Exs \| \vargrad(\xbar_\epoch, z) - \grad \fun (\xbar_\epoch) \|^2 + 3 \Exs \| \vargrad(x, z) - \vargrad(\xbar_\epoch, z) \|^2 \\
&\qquad \qquad \qquad \qquad   + 3 \Exs \| \grad \fun(x) - \grad \fun(\xbar_\epoch) \|^2 \\
& \leq 3 \Exs \| \gradfunbar(\xbar_\epoch) - \fun (\xbar_\epoch) \|^2 + 6  \Lipcon^2 \| x - \xbar_\epoch\|^2 \\
& = \frac{3 \Exs \| \grad \f(\xbar_\epoch, z_1) - \grad \fun (\xbar_\epoch) \|^2}{T}  +  6  \Lipcon^2 \| x - \xbar_\epoch\|^2
\end{align*}
The second last line above follows from the definition of $\vargrad$ and the Lipschitz assumption~\eqref{assn:Lipschitz-gradient},and the last line uses the property uses that $\{ \nabla \f(x, z_i) \}_{i = 1}^\epochsize$ are mutually independent and has the same variance. Finally, we show that as the number of epochs $\epoch$ increase, $\widebar{x}_\epoch$ converges to $\xstar$ at a geometric rate, and as a result, for all $x$ close to $\widebar{x}_\epoch$  
\begin{align*}
\Exs \| \widetilde{\nabla} \f(x, z) - \grad \fun (x) \|^2 \lesssim \frac{\Exs \| \grad \f(\xstar, z_1) - \grad \fun (\xstar) \|^2}{T} \qquad \text{after small number of epochs} \;\; \epoch.
\end{align*}
In other words, after small number of epochs $\epoch$, all points $x$ that are close to $\widebar{x}_\epoch$ are also close to $\xstar$, the variance of gradient estimates $\vargrad(x, z)$ are very similar to the noise variance at the optima $\frac{\Exs \| \grad \f(\xstar, z_1) - \grad \fun (\xstar) \|^2}{T}$ --- the best one could hope for.

\paragraph{Proximal step:}
With this intuition at hand, we are now ready to state our variance reduced proximal-procedure in Algorithm~\ref{Algo:VR-prox}. Our algorithm uses a proximal-step at each step. We use the notation 

\begin{align}
\label{eqn:prox-step}
\prox_{\pen}(x -  \stepsize  \vargrad(x)) = 
\min_{y} \left\lbrace \frac{1}{2} \|y -  \left\lbrace x - \stepsize \vargrad(x) \right\rbrace \|^2 + \pen(y) \right\rbrace  
\end{align}
\begin{algorithm}[ht]
\caption{ $\;\;\;\;\;$ Variance-reduced proximal-gradient (VRPG)
}\label{Algo:VR-prox}
\begin{algorithmic}[1]
\STATE{Given (a) $\nobs$ iid samples $\{z_1, \ldots, z_\nobs\}$ form $\P_0$ (b) Number of epochs $\NumEpoch$, (c) Epoch length
  $\epochLen$ $\;$ (d) Re-centering sample size $\epochsize$}
  (e) Step size $\stepsize > 0$:
  \vspace{10pt} 
  \STATE{Initialize at $\xbar_1$}
\FOR{ epoch $\epiter = 1, 2, \ldots, \NumEpoch$ }
\STATE{Set index $j_\epoch = (\epoch - 1)\cdot (T + K)$.}
\STATE{Use samples $\{ z_i \}_{i = j_\epoch + 1}^{ j_\epoch + \epochsize + \epochLen}$ to construct  $\{\vargrad_{k, \epoch}\}_{k = 1}^K$ as follows
\begin{align*}
\vargrad_{k, \epoch}(x) &:= \grad \f(x, z_{j_\epoch + T + k}) + 
\left\lbrace \gradfunbar(\xbar_\epoch) - \grad \f(\widebar{x}_\epoch, z_{j_\epoch + T + k}) \right\rbrace \\
\text{where} \quad \gradfunbar(x) &= \frac{1}{\epochsize} \sum_{i = 1}^\epochsize \grad \f(x, z_{j_\epoch + i }) 
\end{align*}
 }
\FOR{ k = 1, \ldots K}
\STATE{Set $x_0 = \xbar_\epoch$}
\STATE{Compute proximal update
\begin{align*}
x_{k} = \prox_{\pen}(x_{k - 1} -  \stepsize \vargrad_{k, \epoch}(x_{k - 1}))
\end{align*}
}  
\ENDFOR
\STATE{Set $\xbar_{\epoch + 1} = x_K$ 
  } 
  \ENDFOR \STATE{ Return $x_{\nobs} := \xbar_{\NumEpoch + 1}$ as
  the final estimate}
\end{algorithmic}
\end{algorithm}
%
%

 %%%%%%%%%%%%%%%%%%%%%%%%%%%%%%%%%%%%%%%%%%%%%%%%%%%%%%%%%%%%%%%%%%%%%%%%
 
\subsection{Instance-dependent guarantees}
\label{sec:VRPE-thm-section}
In this section, we state the performance guarantees of the variance reduced proximal method Algorithm~\ref{Algo:VR-prox}.  Given access to $\nobs$ iid samples $\{z_1, \ldots, z_\nobs\}$ from $\P_0$, we 
set 
\begin{align}
\label{eqn:VRPG-tuning}
 \NumEpoch &=  \ceil{\log \nobs}, \quad  \epochsize  = \ceil{\frac{\nobs}{\log \nobs}}, \quad K = \frac{\log(120)}{\log(1/(1 - \muPar^2 / (6 \cdot 8^2 \cdot \Lipcon^2)))} \quad \text{and} \notag \\
 \stepsize &=  \frac{\muPar}{6 \cdot 8^2 \cdot \Lipcon^2} 
\end{align}

\begin{theorem}
\label{thm:VRPG}
Suppose conditions~\eqref{assn:mu-L-funcs}-~\eqref{assn:Lipschitz-gradient} are in force and we have access to $\nobs$ iid samples $\{z_1, \ldots, z_\nobs\}$ from $\P_0$ with $\frac{\nobs}{\log \nobs} \geq 60 \cdot 32 \cdot \frac{\Lipcon^2}{\muPar^2}$. Then, the output $\widehat{x}_{\nobs}$ of Algorithm~\ref{Algo:VR-prox}, obtained with tuning parameters from~\eqref{eqn:VRPG-tuning} and initialization $\widebar{x}_1$, satisfies 
\begin{align}
\label{eqn:min-ell-2-bound}
\Exs\|\widehat{x}_{\nobs} - x^\star \|^2
\leq \frac{\|\widebar{x}_1 - \xstar\|^2}{\nobs^2} + \frac{7\log \nobs}{\nobs} \cdot  \delta^2\left(\frac{\nobs}{\log \nobs}\right)
\end{align}
\end{theorem}
\noindent See Section~\ref{sec:Main-proof} for
the proof of this theorem. 

\vspace{10pt}

\noindent 
A few comments regrading Theorem~\ref{thm:VRPG} are in order.

\paragraph{Non asymptotic nature of upper bound:}
We point out that our upper bound is non-asymptotic in nature, and the bound holds whenever the number of samples $ \nobs \gtrsim \frac{\Lipcon^2}{\muPar^2}$, up to some logarithmic factors. The lower bound $\frac{\Lipcon^2}{\muPar^2}$ is needed to distinguish a $\mu$-strongly function from a constant function. The leading term in our upper bound~\eqref{eqn:min-ell-2-bound} is $\frac{\log \nobs}{\nobs} \cdot \delta^2\left( \frac{\nobs}{\log \nobs } \right)$ --- which differs from our aimed upper bound~\eqref{eqn:UB-banchmark-aim} of $\frac{\delta^2(\nobs)}{\nobs}$ by logarithmic factors. Our bound~\eqref{eqn:min-ell-2-bound} also involves a term which depends on the initial-gap $\|\widebar{x}_1 - \xstar\|^2$; this decencies on the initialization-dependent term decay to zero at a rate $\frac{1}{\nobs^2}$.

\paragraph{Choice of step-sizes:}
Theorem~\ref{thm:VRPG} analyzes a variance reduced algorithm with a constant stepsize $\stepsize$. A similar analysis can be performed with appropriate polynomially decaying stepsizes, i.e, $\stepsize = \frac{c_\omega}{k^\omega}$ for $\omega \in (0, 1)$, and $c_\omega$ is an appropriate universal constant. See for instance the analysis in~\cite{wainwright2019stochastic}. 

\paragraph{Comparison to standard variance reduction methods:}
As we mentioned before, our algorithm is motivated from standard
variance reduction schemes~\cite{johnson2013accelerating,xiao2014proximal} and proceeds in epochs. In usual variance reduction schemes~\cite{johnson2013accelerating,xiao2014proximal}, the number of samples used in epochs increases exponentially with the number of epoch $\epoch$. In our algorithm, however, the number of samples used in each epoch remains same. This change in our algorithm is due to some technical reasons and it contributes to additional factors of $\log(\nobs)$ in our upper bounds. This produces a leading term of  $ \frac{\log \nobs}{\nobs} \cdot \delta^2\left( \frac{\nobs}{\log \nobs } \right)$ instead of $ \delta^2\left( \nobs \right) / \nobs$. It would be interesting to know if these additional logarithmic factors can be removed.  

\paragraph{Removing logarithmic terms in asymptotic limit:} 
It is natural to ask what is the asymptotic behavior of the output of Algorithm~\ref{Algo:VR-prox}. Assuming the output $x_\nobs$ remain bounded, taking number of samples $\nobs \rightarrow \infty$ in the bound~\eqref{eqn:min-ell-2-bound} we see that the upper bound matches the lower bound~\eqref{eqn:minimax-lb} up to a factor of $\log(\nobs)$, and a universal constant. It is interesting to ask whether, these logarithmic factor can be removed, at least asymptotically. Indeed, assuming the iterates $x_\nobs's$ are bounded, from the distributional limit~\eqref{eqn:perturbation-bound} we have that 
\begin{align*}
c_1 \leq \frac{\delta^2(\nobs)}{\nobs} \leq c_2 \qquad \text{for all}
\;\; \nobs \geq \nobs_0. 
\end{align*} 
where $c_1, c_2$ are universal constant.  We can now modify our variance reduced algorithm in the following way. Take $T = \nobs_0$ in epoch $\epoch = 1$, and double the number of samples used in each epoch --- following a standard variance reduction scheme~\cite{xiao2014proximal,johnson2013accelerating}; the choice of $K$ and the stepsize $\stepsize$ remain the same. A simple modification of our current analysis ensures that the output of this modified algorithm achieves the lower bound~\eqref{eqn:minimax-lb} up to universal constants.

\section{Proof of Theorem~\ref{thm:VRPG}}
\label{sec:Main-proof}
To simplify the proof, without loss of generality, we assume that that $\log(\nobs)$, $\nobs/\log(\nobs)$ and the tuning parameters $T$ and $K$ are all positive integers. If this is not the case, the proof can easily be modified by taking appropriate floors, and ceilings. The proof of the theorem involves in three steps.
\begin{itemize}
\item[(a)] Within epoch progress
\item[(b)] Lipschitz property of perturbed problem.
\item[(c)] Between epoch progress
\end{itemize}

\subsection{Within epoch analysis:}
 In Lemma~\ref{lem:SA-lemma}, we characterize the progress of Algorithm~\ref{Algo:VR-prox} within each epoch. 
Within each epoch, we run a constraint stochastic approximation problem. In particular, at epoch $\epoch$ we consider the following perturbed problem  
\begin{align}
\label{eqn:epoch-perturbed-prob-sol} 
  \widehat{x}_\epoch = \arg\min_{x} \left \lbrace  \f(x) +  \left\langle  x, \widebar{\nabla} \f(\xbar_{\epoch}) - \nabla \f(\xbar_\epoch)   \right\rangle  +  \pen(x) \right\rbrace 
  \equiv \arg\min_{x} \left\lbrace \tilde{\f}_\epoch(x) +  \pen(x) \right\rbrace  
\end{align} 
where 
\begin{align*}
\widebar{\nabla} \f(\xbar_{\epoch}) = \frac{1}{\epochsize} \sum_{i = 1}^{\epochsize} \grad \f(\widebar{x}_\epoch, z_{j_\epoch + i}) \qquad \text{with} \quad j_\epoch = (\epoch - 1)*(T + K) 
\end{align*}
In terms of the above notation, $\xbar_{\epoch + 1}$ is the output of a projected stochastic approximation algorithm after $K$ steps. With this set up, we prove  
\begin{lemma}
\label{lem:SA-lemma}
Under assumptions of Theorem~\ref{thm:VRPG} for each epoch $\epoch = 1, \ldots, \NumEpoch$, we have:
\begin{align}
\label{eqn:SA-bound}
\Exs \| \xbar_{\epoch + 1} - \widehat{x}_{\epoch} \|^2_2 \leq \frac{\Exs \| \xbar_{\epoch} - \widehat{x}_{\epoch} \|^2_2}{20}   
\end{align} 
\end{lemma}

\subsubsection*{Proof of Lemma~\ref{lem:SA-lemma}}
Let $x_k$ be the $k^{th}$ iterate of the stochastic approximation within epoch $\epoch$, and we define 
\begin{align*}
g_k = \frac{1}{\stepsize_k} \left( x_{k - 1} - x_{k} \right) 
\end{align*}
For notational convenience we hide the dependence on the epoch $\epoch$, and we use $\widehat{x}$ for $\widehat{x}_{\epoch}$, $\widebar{x}$ for $\widebar{x}_{\epoch}$, and $\tilde{\f}$ for $\tilde{\f}_\epoch$. 
We have 

\begin{align}
\label{eqn:sq-error-expansion}
\|x_{k} - \widehat{x} \|^2 &= \| x_{k - 1} - \widehat{x} \|^2 - 2 \stepsize_k g_k^\top (x_{k - 1} - \widehat{x}) + \stepsize_k^2 \|g_k\|^2 
\end{align}
We next bound the last two terms in the right hand side of the above equation  using the following bound
\begin{align}
\label{eqn:prox-grad-UB}
- g_k^\top (x_{k - 1} - \xhat) + \frac{\stepsize_k}{2} \|g_k\|^2
& \leq  \left \lbrace \tilde{\f}(\xhat) +  \pen(\xhat) - \tilde{\f}(x_k) -  \pen(x_{k})  \right\rbrace - \frac{\muPar}{2} \| x_{k - 1} - \widehat{x} \|^2 - \Delta_k^\top (x_k - \xhat)    
\end{align}
Here, we used the shorthand $\Delta_k$ to denote the noise in the noisy gradient at stage $k$:
\begin{align*} 
\Delta_k = \vargrad_{k, \epoch}(x_{k - 1}) - \nabla \tilde{\f}(x_{k - 1})
\end{align*}
We prove the inequality~\eqref{eqn:prox-grad-UB} shortly. Next note that 
\begin{align}
\tilde{\f}(x_k) +  \pen(x_{k})
&\geq \tilde{\f}(\xhat) +  \pen(\xhat) +  \left\langle \grad \widetilde{\f}(\xhat) +  \partial \pen(\xhat), x_k - \xhat \right\rangle + \frac{\muPar}{2}\|x_k - \xhat \|^2 \notag \\
&= \tilde{\f}(\xhat) +  \pen(\xhat)  + \frac{\muPar}{2}\|x_k - \xhat \|^2  
\label{eqn:func-diff-bound}
\end{align}
The last equality uses $0 \in \grad \widetilde{\f}(\xhat) +  \partial \pen(\xhat)$. 
To bound the term $ -\Delta_k^\top (x_k - \xhat)$, using $\tilde{x}_k = \prox_{\pen}(x_{k - 1} - \stepsize_k \grad \tilde{f}(x_{k - 1}))$,  we have 
\begin{align}
- \Delta_k^\top (x_k - \xhat) &= - \Delta_k^\top (x_k - \tilde{x}_k) - \Delta_k^\top (\tilde{x}_k - \xhat) \notag \\ 
& \leq \| \Delta_k\| \cdot \|x_k - \tilde{x}_k \| -  \Delta_k^\top (\tilde{x}_k - \xhat) \notag \\
&\stackrel{(i)}{\leq} \stepsize_k \|\Delta_k\|^2 -  \Delta_k^\top (\tilde{x}_k - \xhat) \notag \\
&\stackrel{(ii)}{\leq} 8\stepsize_k \Lipcon^2 \|x_{k - 1} - \xhat\|^2 + 2\stepsize_k \|\Delta^\star_k\|^2  -  \Delta_k^\top (\tilde{x}_k - \xhat)  
\label{eqn:inprod-term}
\end{align}
The inequality $(i)$ above uses the contractive property of the proximal operator :
\begin{align*}
\|x_k - \widetilde{x}_k\|_2 &= 
\|\prox_{\pen}(x_{k - 1} - \stepsize_k \vargrad_{k, \epoch}(x_{k - 1})  ) - \prox_{\pen}(x_{k - 1} - \stepsize_k \grad \tilde{f}(x_{k - 1}) ) \| \\
 &\leq \stepsize_k \| \vargrad_{k, \epoch}(x_{k - 1})  - \grad \tilde{f}(x_{k - 1}) \| = \stepsize_k \|\Delta_k\|_2 
\end{align*}
In inequality $(ii)$ we used Cauchy-Schwartz inequality, the $\Lipcon$-Lipschitz property of the maps $\grad \tilde{\f}$ and $\vargrad_{k, \epoch}$ and the following upper bound
\begin{align*}
\|\Delta_k\|_2 &= 
 \| \vargrad_{k, \epoch}(x_{k - 1}) - \nabla \tilde{\f}(x_{k - 1}) \|_2 \\
&\leq \underbrace{\| \vargrad_{k, \epoch}(\xhat) - \nabla \tilde{\f}(\xhat) \|_2}_{:= \|\Delta^\star_k\|_2} + \| \vargrad_{k, \epoch}(\xhat) - \vargrad_{k, \epoch}(x_{k - 1}) \|_2
+ \|\nabla \tilde{\f}(\xhat) - \nabla \tilde{\f}(x_{k - 1})  \|_2 \\
&\leq \|\Delta^\star_k\|_2 + 2 \Lipcon \|x_{k - 1} - \xhat \|_2  
\end{align*}
Substituting bounds~\eqref{eqn:prox-grad-UB},\eqref{eqn:func-diff-bound}, and~\eqref{eqn:inprod-term} in the expression~\eqref{eqn:sq-error-expansion}, and noting $\Exs \left\lbrace \Delta_k^\top (\tilde{x}_k - \xhat) \right\rbrace = 0 $ we have 
\begin{align}
\Exs\|x_{k} - \widehat{x} \|^2 &=  \Exs\| x_{k - 1} - \widehat{x} \|^2 + 2 \stepsize_k \Exs\left\lbrace - g_k^\top (x_{k - 1} - \widehat{x}) + \frac{\stepsize_k}{2} \|g_k\|^2 \right\rbrace \notag \\
&\leq \Exs\| x_{k - 1} - \widehat{x} \|^2 + 2 \stepsize_k \left\lbrace -  \frac{\muPar}{2}  \Exs\|x_{k - 1} - \xhat\|^2 - \frac{\muPar}{2} \Exs\|x_k - \xhat \|^2 + 2\stepsize_k \Exs\|\Delta_k^\star\|^2 + 8 \stepsize_k \Lipcon^2 \Exs\| x_{k - 1} - \xhat \|_2^2 \right\rbrace 
\label{eqn:sq-error-bound-2}
\end{align}
Taking $\stepsize_k \leq \frac{\muPar}{6 \cdot 8^2\Lipcon^2}$ we have 
\begin{align*}
16 \stepsize_k^2 \Lipcon^2 - \muPar \stepsize_k \leq 0, \qquad \text{and} \qquad  \frac{1}{1 + \muPar \stepsize_k} \leq \left(1 - \muPar \stepsize_k \right) = \left(1 - \frac{\muPar^2}{6 \cdot 8^2 \Lipcon^2}\right) 
\end{align*}  
We have 
\begin{align*}
\Exs\|x_{k} - \widehat{x} \|^2 &\leq  
\frac{1}{1 + \muPar \stepsize_k} 
\Exs \|x_{k - 1} - \widehat{x} \|^2 
+ 4 \stepsize_k^2 \Exs \| \Delta_k^\star \|^2   
\end{align*}
Finally, using the Lipschitz property of $\grad \f(x,z)$ and $\grad \f(x)$ from Assumptions~\eqref{assn:Lipschitz-gradient} and~\eqref{assn:mu-L-funcs}, respectively, we have 
\begin{align*}
\| \Delta_k^\star \|_2
= \| \nabla \f(\xhat, z_{j_\epoch + T + k}) - 
\nabla \f(\widebar{x}, z_{j_\epoch + T + k}) + \nabla \f(\widebar{x}) - \nabla \f(\xhat)  \|_2  \leq 2 \Lipcon \| \xhat - \widebar{x} \|_2 
\end{align*}
Putting together the pieces we obtain at the following recursion
\begin{align}
\Exs\|x_{k} - \widehat{x} \|^2 &\leq \left(1 - \frac{\muPar^2}{6 \cdot 8^2 \cdot \Lipcon^2} \right) 
\Exs \|x_{k - 1} - \widehat{x} \|^2 
+   \cdot \frac{16 \cdot\muPar^2}{6^2  \cdot 8^4 \cdot \Lipcon^2} \| \xhat - \widebar{x} \|_2^2 
\end{align}
Recursing the last equation $k$-times with $k \geq \frac{\log(120)}{\log(1/(1 - \muPar^2 / 6 \cdot 8^2 \cdot \Lipcon^2))}$ we have 
\begin{align*}
\Exs\|x_{k} - \widehat{x} \|^2 &\leq \left(1 - \frac{\muPar^2}{6 \cdot 8^2 \cdot \Lipcon^2} \right)^k 
\Exs \|x_0 - \widehat{x} \|^2 
+  \frac{16}{6 \cdot 8^2} \cdot \frac{\muPar^2 / (6 \cdot 8^2 \cdot \Lipcon^2)}{\muPar^2 / (6 \cdot 8^2 \cdot \Lipcon^2)} \Exs \| \xhat - \widebar{x} \|_2^2 \notag \\
&=  \left(1 - \frac{\muPar^2}{6 \cdot 8^2 \Lipcon^2} \right)^{k} 
 \Exs \|\widebar{x} - \widehat{x} \|^2 + \frac{1}{24} \Exs \|\widebar{x} - \widehat{x} \|^2 \\
 &\leq \frac{1}{120} \Exs \|\widebar{x} - \widehat{x} \|^2 + \frac{1}{24} \Exs \|\widebar{x} - \widehat{x} \|^2 \leq \frac{\Exs \|\widebar{x} - \widehat{x} \|^2}{20}
\end{align*}
Here, we used $\left(1 - \frac{\muPar^2}{6 \cdot 8^2 \Lipcon^2} \right)^k \leq 1/120$ for $k \geq \frac{\log(120)}{\log(1/(1 - \muPar^2 / 6 \cdot 8^2 \cdot  \Lipcon^2))}$. It remains to prove the bound~\eqref{eqn:prox-grad-UB}. 

\subsubsection*{Proof of bound~\eqref{eqn:prox-grad-UB}:}
The bound follows from the following \emph{key} lemma on the progress of proximal operators. 
\begin{lemma}
\label{lem:proximal-progress-lemma}
Given a function For any stepsize $0< \stepsize \leq \frac{1}{\Lipcon}$ we have 
\begin{align}
\f(y) + \pen(y) \geq f(x^+) + \pen(x^+) +  g^\top (y - x) + \frac{\stepsize}{2} \|g\|^2 + 
\frac{\muPar}{2}\|y - x\|^2 + \Delta^\top(x^+ - y)  
\end{align}
where $x^+ = \prox_{\pen}(x - \stepsize v)$ and $\Delta = v - \grad \f(x)$. 
\end{lemma}
Given this lemma at hand, we first note that the perturbed function $\ftilde$ is 
$\muPar$-strongly convex and $\Lipcon$-smooth, and the stepsize $\stepsize = \frac{\mu}{16\Lipcon} \cdot \frac{1}{\Lipcon} \leq 1/\Lipcon$. Now applying Lemma~\ref{lem:proximal-progress-lemma} 
on the function $\ftilde$, $y = \widehat{x}$, $x = x_{k - 1}$, $\stepsize = \stepsize_k$ and $v = \vargrad_{k, \epoch}(x_{k - 1})$ yields the bound~\eqref{eqn:prox-grad-UB}. We provide a proof of Lemma~\ref{lem:proximal-progress-lemma} in Appendix~\ref{lem:proximal-progress-proof}.

\subsection{Lipschitz property of perturbed problem solutions:}
The second ingredient in the proof is a Lipschitz property of the perturbed problem. Consider two following perturbed problems
\begin{align*}    
x_{\tilde{\f}} = \arg\min_{x} \left\lbrace \f(x) + \left\langle  x, \widebar{\nabla} \f(\xbar) - \nabla \f(\xbar)   \right\rangle + \pen(x) \right\rbrace \equiv \arg\min_x \tilde{\f}(x)  \qquad \text{and} 
\end{align*}

\begin{align*}    
x_{\tilde{g}} = \arg\min_{x} \left\lbrace \f(x) 
+ \left\langle  x, \widebar{\nabla} \f(x^\star) - \nabla \f(x^\star) \right\rangle + \pen(x) \right\rbrace \equiv \arg\min_x \gtilde(x) 
\end{align*}
% %

\noindent  Here, $\widebar{x}$ is any point independent of the randomness in $\widebar{\nabla}\f$, and $x^\star$ is the (unique) solution of the problem~\eqref{eqn:constraint-opt}.
The $\muPar$-strong convexity of the function $\f$ ensures that the function $\ftilde$ satisfies the following two  growth condition
\begin{align*}
\ftilde(x) &\geq \ftilde(x_{\ftilde}) + \frac{\muPar}{2}\|x -  x_{\ftilde} \|^2 \qquad \text{for all} \;\; x.  
\end{align*} 
We have 

\begin{align*}
\frac{\muPar}{2}\|x_{\gtilde} - x_{\ftilde}\|^2 \leq  \ftilde(x_{\gtilde}) - \ftilde(x_{\ftilde}) &=    \left\lbrace \ftilde(x_{\gtilde}) - \gtilde(x_{\gtilde}) \right\rbrace + \underbrace{\left\lbrace \gtilde(x_{\gtilde}) - \gtilde(x_{\ftilde}) \right\rbrace}_{\leq 0} - \left\lbrace \ftilde(x_{\ftilde}) - \gtilde(x_{\ftilde}) \right\rbrace \\
% %
&\leq \left\lbrace \ftilde(x_{\gtilde}) - \gtilde(x_{\gtilde}) \right\rbrace  - \left\lbrace \ftilde(x_{\ftilde}) - \gtilde(x_{\ftilde}) \right\rbrace \\
&= \left\langle x_{\gtilde} - x_{\ftilde}, [\widebar{\nabla} \f(\widebar{x}) - \widebar{\nabla} \f(x^\star)]
-[
\nabla \f(\widebar{x}) - \nabla \f(x^\star)] \right\rangle \\
&\leq \|x_{\gtilde} - x_{\ftilde} \| \cdot \| [\widebar{\nabla} \f(\widebar{x}) - \widebar{\nabla} \f(x^\star)]
-[
\nabla \f(\widebar{x}) - \nabla \f(x^\star)] \|
\end{align*}
Rearranging, last bound and taking expectation we have 

\begin{align}
\label{eqn:Lipschitz-property-of-sols}
\Exs \|x_{\gtilde} - x_{\ftilde}\|^2
\leq \frac{4 \Exs \| [\widebar{\nabla} \f(\widebar{x}) - \widebar{\nabla} \f(x^\star)]
-[
\nabla \f(\widebar{x}) - \nabla \f(x^\star)] \|^2}{\muPar^2} \leq \frac{16 \Lipcon^2 \|\widebar{x} - x^\star\|^2}{T\muPar^2}
\end{align}
The last inequality above follows from the $\Lipcon$-Lipschitz property of the gradient $\grad \f$.

\subsection{Between epoch progress:}
We complete the proof by quantifying the progress over epochs. Recalling the definition of $\widehat{x}_\epoch$ and $x_{\nobs/\log \nobs}^\star$ from~\eqref{eqn:epoch-perturbed-prob-sol} and~\eqref{eqn:tilded-penalized-version}, respectively, we deduce the following bound 
\begin{align}
\Exs\|\xbar_{\epoch + 1} - \xstar\|^2
& \leq 2 \Exs\|\xbar_{\epoch + 1} - \widehat{x}_\epoch \|^2 + 2 \Exs\|\widehat{x}_\epoch - \xstar \|^2 \notag \\
& \stackrel{(i)}{\leq}  \frac{2\Exs \|\xbar_{\epoch} - \widehat{x}_\epoch \|^2 }{20} +  2\Exs \|\widehat{x}_\epoch - \xstar \|^2 \notag \\ 
&\leq \frac{\Exs \|\xbar_{\epoch} - \xstar \|^2}{5} + \frac{\Exs \|\widehat{x}_\epoch - \xstar \|^2}{5} + 2\Exs \|\widehat{x}_\epoch - \xstar \|^2 \notag \\
&= \frac{\Exs \|\xbar_{\epoch} - \xstar \|^2}{5} + \frac{11}{5} \cdot \Exs \|\widehat{x}_\epoch - \xstar \|^2
\label{eqn:bound-1}
\end{align}
Inequality $(i)$ follows from Lemma~\ref{lem:SA-lemma}. Next, we have 
\begin{align}
 \Exs \|\widehat{x}_\epoch - \xstar \|^2 &\leq 2\Exs \|\widehat{x}_\epoch - x_{\nobs/\log \nobs}^\star \|^2  + 2 \Exs \|x_{\nobs/\log \nobs}^\star - \xstar \|^2 \notag \\
 &\leq \frac{2\cdot 16 \Lipcon^2}{T \muPar^2} \cdot  \Exs \|\xbar_{\epoch} - \xstar \|^2  + 2 \Exs \|x_{\nobs/\log \nobs}^\star - \xstar \|^2 \notag \\
 &\leq \frac{5\cdot\Exs \|\xbar_{\epoch} - \xstar \|^2}{11\cdot20} + 2 \Exs \|x_{\nobs/\log \nobs}^\star - \xstar \|^2 \label{eqn:bound-2}
\end{align}
The second last line above used the bound~\eqref{eqn:Lipschitz-property-of-sols} and noting that $T = \frac{\nobs}{\log \nobs}$. The last line uses the fact $T = \frac{\nobs}{\log \nobs} \geq \frac{60\cdot 32 \cdot \Lipcon^2}{\muPar^2}$. Combining bounds~\eqref{eqn:bound-1} and~\eqref{eqn:bound-2}and recursing the bound for $\NumEpoch$ we obtain
\begin{align*}
\Exs\|\widebar{x}_{\NumEpoch + 1} - \xstar\|^2 &\leq \frac{\Exs\|\widebar{x}_{\NumEpoch} - \xstar\|^2}{4} + 5 \cdot \Exs \|x_{\nobs/\log \nobs}^\star - \xstar \|^2 \\ 
&\leq \frac{\Exs\|\widebar{x}_{1} - \xstar\|^2}{4^\NumEpoch} + 7 \cdot \Exs \|x_{\nobs/\log \nobs}^\star - \xstar \|^2 \\
& = \frac{\Exs\|\widebar{x}_{1} - \xstar\|^2}{\nobs^2} + \frac{7 \log(\nobs)}{\nobs} \cdot \delta^2\left( \frac{\nobs}{\log \nobs}\right).  
\end{align*}
This completes the proof of Theorem~\ref{thm:VRPG}.

\section{Discussion}
In this paper we study the non-asymptotic performance of a natural a variance-reduced algorithm (VRPG) for a stochastic approximation problem under constraints. An attractive feature of  our results are that they are  instance-dependent in nature. Meaning our guarantees reflect the complexity of the loss function, the variability of the noise, the geometry constraint set. Additionally, we show that when the number of samples $\nobs \rightarrow \infty$, our upper bound matches the celebrated H\'a jek-Le Cam local minimax lower bound up to a factor of $\log (\nobs)$ and some universal constants.

Our work leaves open a number of open questions. For instance, there is a significant body of work which construct a non-asymptotic local minimax lower bound by comparing the given problem with its hardest alternative in a local shrinking  neighborhood of the given problem~\cite{cai2015framework}. Indeed, in absence of the constraint set (i.e. $\Xset = \real^\dims$), variance reduced schemes --- the ones similar to the one studied in our paper --- achieve such non-asymptotic local minimax lower bound~\cite{khamaru2021instance,khamaru2021instanceQ}. It would be interesting to see 
if such results can be extended to our setting where $\Xset \neq \real^\dims$.

\subsection*{Acknowledgments}
This research was partially supported by National Science Foundation Grant
DMS-2311304 to KK. The author thank Feng Ruan and Qiyang Han for helpful discussions on the paper.

\bibliographystyle{abbrv} \bibliography{ref}

\appendix

\section{Technical lemma}
In this section we prove a technical lemma that were used in the main section of the proof.  

\subsection{Proof of Lemma~\ref{lem:proximal-progress-lemma}:}
\label{lem:proximal-progress-proof}
This Lemma is directly borrowed from~\cite[Lemma 3]{xiao2014proximal}. We provide a proof of this result for completeness. Recall 
that
\begin{align*}
x^+ = \arg\min_y \left\lbrace \frac{1}{2}\|y - (x - \stepsize v)\|^2 + \pen(y)  \right\rbrace 
\end{align*}
Thus, there exists $u \in \partial \pen(x^+)$ such that 
\begin{align}
\label{eqn:prox-stat-eqn}
x^+ = (x - \stepsize v) - \stepsize \cdot u
\end{align}
Now using $\mu$-strong convexity of $\f$ and convexity of the penalty function $\pen$ we have 
\begin{align*}
\f(y) + \pen(y) 
\geq \f(x) + \grad \f(x)^\top (y - x) + \frac{\mu}{2}\|y - x\|^2
+ \pen(x^+) + u^\top(y - x^+)  
\end{align*} 
Additionally, using $\Lipcon$ smoothness of the function $\f$, and the expression $x - x^+= \stepsize g$ we have
\begin{align}
\f(y) + \pen(y) 
&\geq \f(x^+) +  \pen(x^+) - \frac{L\stepsize^2}{2}\|x^+ - x\|^2 + \frac{\mu}{2}\|y - x\|^2 \notag \\
& \qquad - \grad \f(x^+)(x^+ - x) +\grad \f(x)^\top (y - x) 
   + u^\top(y - x^+) \notag \\
& = \f(x^+) +  \pen(x^+) - \frac{L\stepsize^2}{2}\|x^+ - x\|^2 + \frac{\mu}{2}\|y - x\|^2 \notag \\
& \qquad  \qquad +g^\top(y - x) + \stepsize \|g\|^2 + \Delta^\top(x^+ - y) 
\label{eqn:descent-lemma-step}
\end{align}
The last line above uses the following simplification using the relations $\Delta  = v - \grad \f(x)$, $u = g - v$, and $x - x^+ = \stepsize g$ 
\begin{align*}
&- \grad \f(x^+)(x^+ - x) +\grad \f(x)^\top (y - x) 
   + u^\top(y - x^+) \\
&=g^\top(y - x^+) + (v - \grad \f(x))^\top(x^+ - y) \\  
&=g^\top(y - x) + \stepsize \|g\|^2 + \Delta^\top (x^+ - y)     
\end{align*}
Finally, using the stepsize condition $0 < \stepsize \leq \frac{1}{\Lipcon}$ in the bound~\eqref{eqn:descent-lemma-step} we conclude 
\begin{align*}
\f(y) + \pen(y) \geq \f(x^+) + \pen(x^+) + \frac{\stepsize}{2}\|g\|^2
+ \frac{\mu}{2}\|y - x^+\|^2 + \Delta^\top(x^+ - y).
\end{align*}
This completes the proof of Lemma~\ref{lem:proximal-progress-lemma}.

\end{document}